\newcommand{\R}{\mathbb{R}} 
 
\newcommand{\Z}{\mathbb{Z}}

\mathchardef\varepsilon="010F
\mathchardef\epsilon="0122
\mathchardef\vartheta="0112
\mathchardef\theta="0123
\mathchardef\varrho="011A
\mathchardef\rho="0125
\mathchardef\varphi="011E     
\mathchardef\phi="0127
\renewcommand \emptyset \varnothing
\documentclass[12pt]{article}            
\usepackage[all]{xy}
\usepackage[latin1]{inputenc}        
\usepackage[dvips]{graphics}
\usepackage[dvips]{graphicx}
\usepackage{amsfonts}
\usepackage{amssymb}
\usepackage{amsmath}
\usepackage{amstext}
\usepackage{amsbsy}
\usepackage{amsopn}
\usepackage{amsthm}
\usepackage{amscd}
\usepackage{amsxtra}
\usepackage{mathrsfs}
\usepackage{upref}
\usepackage{fullpage}
\title{On geodesic envelopes}
\author{Gianmarco Capitanio}

\begin{document}        

\theoremstyle{plain}
\newtheorem*{theorem*}{\bf Theorem}
\newtheorem{theorem}{\bf Theorem}
\newtheorem*{conjecture*}{\bf Conjecture}
\newtheorem{lemma}{\bf Lemma}
\newtheorem*{lemma*}{\bf Lemma}
\newtheorem{proposition}{\bf Proposition}
\newtheorem*{proposition*}{\bf Proposition}
\newtheorem{corollary}{\bf Corollary}
\newtheorem*{corollary*}{\bf Corollary}
\theoremstyle{definition}
\newtheorem*{definition*}{\bf Definition}
\newtheorem*{definitions*}{\bf Definitions}

\newtheorem{example}{\bf Example}
\newtheorem*{example*}{\bf Example}
\theoremstyle{remark}
\newtheorem*{remark*}{\bf Remark}
\newtheorem{remark}{\bf Remark}
\newtheorem*{remarks*}{\bf Remarks}
\newtheorem{remarks}{\bf Remarks}

\maketitle

\begin{abstract}
We give a global description of envelopes of geodesic 
tangents of regular curves in (not necessarily convex) Riemannian
surfaces.  
We prove that such an envelope is the union of the curve itself,  
its inflectional geodesics and its tangential caustics (formed by the conjugate
points to those of the initial curve along the tangent geodesics). 
Stable singularities of tangential caustics and geodesic envelopes are discussed. 
We also prove the (global) stability of tangential caustics of
close curves in convex closed surfaces under small deformations of the initial
curve and of the ambient metric.  
\end{abstract}
 

\section{Introduction}

The {\it geodesic envelope} of a regular curve in a (complete smooth)
Riemannian surface is the envelope of its geodesic  tangents.  
In this note we describe geodesic envelopes from a global viewpoint. 
For this, we consider the {\it tangential caustics} of the curve, formed by 
the conjugate points (of different order) to those of the curve along its tangent
geodesics.  

Recall that two points are conjugate if the second point is the
intersection of two arbitrary close geodesics issuing from the first
point; the caustics associated to a point are the sets of the $n$-th
conjugate points to the initial point (the ordering of the conjugate
points along a geodesic being well-defined, due 
to the positiveness of the injectivity radius). 

The tangential caustic is not the usal caustic defined in
Optics, which is the envelope of the 
normal geodesics of a curve. 
It is a different generalisation of the usual caustic of a point (formed by its
conjugate points), introduced by Poincaré in \cite{poincare}. 

Our first result is the following global description of geodesic envelopes.

\begin{theorem}\label{thm:1}
The geodesic envelope of any regular curve in a complete Riemannian surface 
is the union of the curve, its inflectional geodesics and its
tangential caustics.  
\end{theorem} 

This result generalizes to Riemannian manifolds some of Thom's
results about envelopes of $1$-parameter families of lines in 
the projective plane (see \cite{thom}).   
We also prove that the only generic and stable singularities of
tangential caustics are semicubic cusps and transversal
self-intersections; for geodesic envelopes, we have to add local
second-order self-tangencies (which are not stable for general
envelopes).   

The above theorem and the classification of stable singularities of
tangential caustics lead to the following stability statement. 

\begin{theorem}\label{thm:2}
Each $p$-tangential caustic of a regular closed curve in a strictly convex
closed Riemannian surface is generically stable under
small enough deformations of the initial curve and of the ambient
metric (the required smallness depending on the order of the
caustic).  
\end{theorem}

The stability means that there exists a diffeomorphism of the ambient
surface transforming the unperturbed $n$-caustic into the perturbed
one. 
``Generically'' means (here and throughout the paper) 
that the curves for which the statement does not
hold form a residual set for Whitney topology in the set of all the
curves in the surface (see \cite{ST1}, §3).  

The framework here is the same as several classical subjects in Geometry,
Optics and Calculus of Variations, going back for instance to
Archimedes, Huygens, Barrow and Jacobi, as the study of evolutes
of curves and caustics of ray systems. 
The relation between singularities of ray systems, their
caustics, wave fronts, Legendre transformations and reflection groups
was discovered by Arnold (see \cite{arnold1}), in the setting of
Symplectic and Contact Geometry, and further developed by
O.V. Lyashko, A.B. Givental, O.P. Shcherbak (see \cite{ST2}). 

\medskip

\noindent {\bf Acknowledgements.} Research funded by INdAM.

\section{Preliminary definitions}\label{sct:2}

Let $(M,g)$ be a complete smooth Riemannian surface with metric $g$ and let 
$\gamma:I\rightarrow (M,g)$ be a regular curve parameterized by
arc-length, where $I$ is an open interval of the circle $S^1$.   
We endow as usual the manifold of all such maps with Whitney topology. 
For each $\xi\in I$, denote by $\Gamma_\xi:\R\rightarrow M$ the
unit speed geodesic, tangent to
$\gamma$ at $\gamma(\xi)$ (and oriented as $\gamma$ at this point). 

We say that $\gamma$ is {\it convex} at $\xi$ if $\ddot \gamma(\xi)\not=0$;
in this case, the tangency order of $\gamma$ with the tangent geodesic 
$\Gamma_\xi$ at $\gamma(\xi)$ is $1$. 
If $\ddot \gamma(\xi)=0$, then $\xi$ is said to be an {\it inflection} of $\gamma$.
The inflection is {\it simple} when the curve has a second-order
tangency with the inflectional geodesic. 

The {\it graph} of the family of geodesics $\{\Gamma_\xi :\xi\in I\}$
is the immersed surface  
$$\Phi:= \left\{\big(\Gamma_\xi(t),\xi\big) : \xi\in I, t\in\R\right\}\subset
M\times I\ .$$

\begin{definition*} 
The {\it geodesic envelope} of $\gamma$ is the apparent contour of 
the graph in $M$, that is, the critical value set of the graph's
projection $\Phi\rightarrow M$ defined by $(P,\xi)\mapsto P$.
\end{definition*}
 
Note that the image of $\gamma$ is contained in the geodesic envelope.

Consider an embedded local component $\bar \delta\subset \Phi$ of the critical
set of the graph's projection; denote by $\delta$ its
projection in $M$.   
Let $(P,\xi)\in \bar\delta$. 
The branch $\delta$ is said to be a {\it geometric envelope} at
$P=\Gamma_\xi(t_0)$  if the curve $t\mapsto (\Gamma_\xi(t),\xi)$ 
intersects $\bar\delta$ only at $t=t_0$ (in a small neighbourhood of
$t_0$).  

\begin{definition*}
The {\it tangential caustic} of a regular curve $\gamma$ in a complete
Riemannian surface $(M,g)$ is the set $\Sigma$ 
formed by the conjugate points to those of $\gamma$ along its tangent
geodesics.  
\end{definition*}

As the caustic of a point, the tangential caustic splits into
components $\Sigma_p$, formed by the $p$-th conjugate points.  
Here $p$ belongs to $\Z$, the conjugate points of negative
order lying on the tangent geodesic in the direction opposed to that
of $\gamma$ at the tangency point. 
The curve $\gamma$ can be considered as the $0$-tangential caustic.

\begin{example}
A curve in a surface with non-positive curvature (as for instance
the Euclidian plane or the Lobachevsky plane)  has no tangential
caustics of order $p\not=0$ (see e.g. \cite{docarmo}). 
In this case, Theorem \ref{thm:1} says that the envelope of the
geodesic tangents of a curve is the union of the curve with its
inflectional tangents.   
\end{example}

\begin{example}
Let $(M,g)$ be a Riemannian surface whose curvature is everywhere bounded by two given positive constants. 
Then, every regular curve in $M$ has infinitely many tangential caustics (see \cite{docarmo}); if the curve is closed, then every tangential caustic is also closed.
For instance, the $p$-th tangential caustic of a curve in the standard sphere $S^2$ is the curve itself for even $p$ and the antipodal curve otherwise.  
\end{example}

For $p\in \Z$, let $\tau_p(\xi)$ be the distance along $\Gamma_\xi$
between $\gamma(\xi)$ and its $p$-th conjugated point (provided that
it exists).
In particular, $\tau_0\equiv 0$.     
Each map $\tau_p$ is continuous, and defined on a subset $I_p$ of
$I$, possibly empty. 
We have $I_{p+1}\subseteq I_{p}$ for $p\geq0$, $I_{p}\subseteq I_{p+1}$ otherwise.

\begin{lemma*}
The map $\phi_p:I_p\rightarrow M$, sending $\xi$ to  $\Gamma_\xi(\tau_p(\xi))$,   
is a smooth parameterization of $\Sigma_p$ (which is hence a smooth curve). 
\end{lemma*}

\begin{proof}
Let $P$ and $Q$ be two consecutive conjugated points along a geodesic $c$. 
Then, there exists a non-trivial Jacobi field along $c$, orthogonal at
every point to the geodesic and vanishing at $P$ and $Q$ (see \cite{docarmo}).  
Consider a smooth deformations $P(\lambda)$ and $c(\lambda)$ of $P$
and $c$: $c(\lambda)$ is a geodesic issuing from $P(\lambda)$ for
every small enough $\lambda$.  
Due to the regularity Theorem about solutions of 
second order differential equations, depending on parameters, we 
obtain that, for every $\lambda$ small enough, there exists a
non-trivial orthogonal Jacobi field along $c(\lambda)$, vanishing to
$P(\lambda)$ and to some other point $Q(\lambda)$.  
This point is therefore the conjugate point
of $P(\lambda)$  along the geodesic $c(\lambda)$, and it depends
smoothly on the parameter $\lambda$, provided that the deformation is
small enough. 
\end{proof}

\section{Proof of Theorem \ref{thm:1}}\label{sct:3}
The proof of Theorem \ref{thm:1} is subdivided into three steps. 
In the first two steps we show that the inflectional tangents
and the tangential caustics of a curve are contained in its geodesic
envelope.  
Finally we show that this envelope contains nothing else. 

\begin{proposition}\label{prop:1} 
The geodesic envelope of a regular curve in a Riemannian surface 
contains its inflectional geodesics (as non geometric components).
\end{proposition}

\begin{proof}
Let $P=\gamma(\xi_0)$ be a simple inflection of $\gamma$ and 
let $Q=\Gamma_{\xi_0}(T)$ be the first conjugate point to $P$ along 
$\Gamma_{\xi_0}$ (the argument can be easily adapted to the case where $P$ has no
conjugated points). 
We can suppose $T>0$, and change the orientation of $\gamma$ to
consider the conjugate points in the other direction.

Fix two geodesic balls $B_P$ and $B_Q$, centered at $P$ and $Q$
respectively, of radius $r>0$ arbitrary small.  
For $\xi\rightarrow \xi_0$ with $\xi\not=\xi_0$, the geodesics $\Gamma_\xi$
meet the geodesic $\Gamma_{\xi_0}$ in the interior of $B_P$ and then meet
the ball's frontier $\partial B_P$ at $\Gamma_\xi(t_1(\xi))$, where
$t_1(\xi)=r+o(\xi-\xi_0)$; all these intersections belong to the 
closure of the same connected component of 
$\partial B_P\smallsetminus\{\Gamma_{\xi_0}(\pm r)\}$.
Next, the geodesics $\Gamma_\xi$ intersect $\partial B_Q$ at 
$\Gamma_\xi(t_2(\xi))$, where $t_2(\xi)=T-r+o(\xi-\xi_0)$, 
and $\Gamma_{\xi_0}$ in $B_Q$. 

Hence, the geodesic segments
$\Gamma_\xi(]t_1(\xi),t_2(\xi)[)$ do not intersect 
$\Gamma_{\xi_0}(]r,T-r[)$. 
Thus, their lifting in the graph $\Phi$ form a smooth fold-like surface, 
whose apparent contour in $M$ 
is $\Gamma_{\xi_0}(]r,T-r[)$ for every $r>0$ arbitrary small, so 
$\Gamma_{\xi_0}([0,T])$ is contained in the geodesic envelope.  
The same argument shows that the envelope contains the 
$\Gamma_{\xi_0}$ geodesic segment between any two consecutive
conjugate points and, therefore, the whole geodesic $\Gamma_{\xi_0}(\R)$. 

Finally, if $\gamma$ has non-simple inflections, we consider 
a deformation $\gamma^\lambda$ of the curve $\gamma$, having only simple 
inflections for $\lambda\not=0$. 
The statement holding for the perturbed envelopes, by continuity we
get the claim for the unperturbed envelope.  
This ends the proof, since the non geometricity
of the inflectional geodesic is clear. 
\end{proof}

\begin{proposition}\label{prop:2}
The tangential caustic of a regular curve is contained in its
geodesic envelope.
\end{proposition}

\begin{proof} 
Fix $p\in\Z$ such that $I_p$ is not empty. 
Let $\xi_0\in I$, such that $\tau_p(\xi_0)\in I_p$. 
By the very definition of conjugate point, for every $\epsilon >0$ there exists a $p$-depending constant $C$ such
that for every $|\xi-\xi_0|<C$ the geodesics $\Gamma_{\xi_0}$ and
$\Gamma_\xi$ intersect each other at some points
$\Gamma_{\xi_0}(t_i)$, $i=0,\dots,p$, where
$|t_i-\tau_i(\xi_0)|<\epsilon$. 
By continuity, $t_i\rightarrow \tau_i(\xi_0)$ for
$\xi\rightarrow\xi_0$. 
Hence, every point of the $p$-th tangential caustic
$\Sigma_p$ is the intersection of infinitesimally close curves of the
family of tangent geodesics. 
Therefore, $\Sigma_p$ is contained in the geodesic envelope. 
\end{proof}

\begin{remarks*}
(1) Let us recall that there exists a different definition of envelope: a
point belongs to the ``naif envelope'' of a family of curves if it is
the intersection of infinitesimally close curves of the
family. 
In the preceding proof, we have actually shown that $\Sigma_p$ is
contained in the ``naif envelope'' of the geodesic tangents of
$\gamma$. 
But it is well known (see e.g. \cite{BG}, §5) that the ``naif
envelope'' is contained in Thom's envelope (being sometimes
different).  

\noindent (2) Each tangential caustic is a
geometric branch of the geodesic envelope.

\noindent (3)
The caustic of a point can be viewed as the envelope of the geodesics
issuing from it. 
A similar characterization holds for the tangential caustics of a curve.  
Indeed, $\Sigma_p\cup \Sigma_{-p}$ is the envelope of the pencil,
parameterized by $\xi$, of the usual $|p|$-caustics of the points
$\gamma(\xi)$. 
\end{remarks*}

The proof of Theorem \ref{thm:1} is completed by the following fact. 

\begin{proposition}
Every point of the geodesic envelope of a regular curve $\gamma$
belongs to its tangential caustics or to an inflectional geodesic
tangent. 
\end{proposition}
 
\begin{proof}
Fix a point $\Gamma_\xi(t)$ of the geodesic envelope of $\gamma$.  
Thus, $(\Gamma_\xi(t),\xi)$ is a
critical point of the graph's projection $\pi$.  
Suppose that this point belongs to a regular geometric branch of the
envelope (if the branch is not regular, we can make it regular by an
arbitrary small deformation of $\gamma$).  
In this case the geodesics $\Gamma_\xi$ and
$\Gamma_{\xi+\epsilon}$ intersect each other at two points
$\Gamma_\xi(o(\epsilon))$ and $\Gamma_\xi(t+o(\epsilon))$ for
$\epsilon$ arbitrary small.  
Passing to the limit for $\epsilon\rightarrow 0$ we obtain either that
$\Gamma_\xi(0)$ and $\Gamma_\xi(t)$ are conjugated or that $t=0$.  
In both cases $\Gamma_\xi(t)$ belongs to a tangential caustic of
$\gamma$.   

On the other hand, if the envelope is not geometric at
$(\Gamma_\xi(t),\xi)$, then it is clear that the whole geodesic
$\Gamma_\xi$ is contained in the geodesic envelope of $\gamma$.  
As  in the proof of Proposition \ref{prop:1}, this is possible if and
only if $\gamma$ has an inflection at $\xi$.  
\end{proof}
 
\section{Proof of Theorem \ref{thm:2}}
We start describing the singularities of tangental caustics.

\begin{proposition}\label{prop:4}
Each tangential caustic $\Sigma_p$ of a regular curve $\gamma$ 
is regular outside a subset $X_p\subset I_p$ of $I_p$, which is
generically discrete.  
The generic singularities of the tangential caustics are semicubic
cusps and transversal self-intersections; these singularities are
stable under small enough deformations of the curve $\gamma$ and of
the ambient metric $g$.  
\end{proposition} 

\begin{figure}[h]
  \centering
   \scalebox{.7}{\input{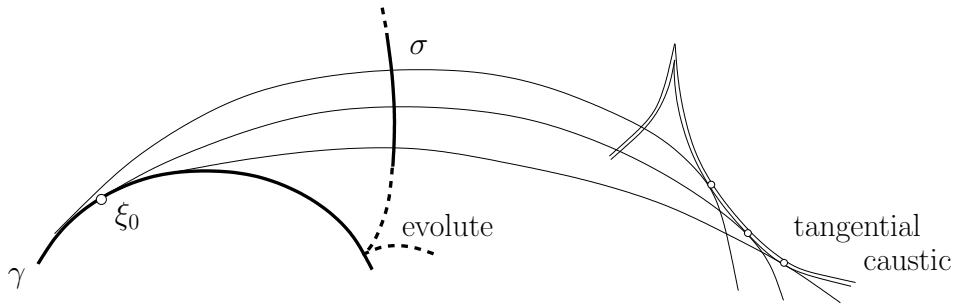}}  
\caption{The tangent geodesics of a curve as a system of rays.}
 \label{fig:2}
\end{figure}

\begin{proof}
Recall that a {\it system of rays} emanating from a given embedded
curve is the pencil of the normal geodesics of the curve. 
The envelope of these normal geodesics is called the {\it caustic of
the system of rays}.   

Fix a point $\xi_0\in I$, such that $\ddot \gamma(\xi_0)$ is not
vanishing. 
Then there exists some evolvent $\sigma$ of $\gamma$, such that the
tangent geodesics to $\gamma$, $\Gamma_\xi$, form a system of rays
emanating from $\sigma$, provided that $\xi$ is close enough to
$\xi_0$ (see fig. \ref{fig:2}). 
The generic singularities of the caustics of systems of rays are
semicubic cusps and transversal self-intersections (see \cite{ST2}). 
These singularities are stable under small deformations of the system
of rays.  
Thus, each curve $\phi_p$ is regular outside a generically discrete
subset of $I_p$.  
\end{proof}

The genericity required in the Lemma provides the
genericity of the Lagrangian surface associated to the system of rays,
see \cite{ST2} and therein references.   

\begin{proposition}\label{prop:5}
Let $\gamma(\xi_0)$ be a simple inflection of $\gamma$. 
Then, for every fixed $p\in \Z$, $\Sigma_p$ has generically a simple
inflection at the corresponding point $\phi_p(\xi_0)$. 
\end{proposition}

\begin{proof}
Generically $\Sigma_p$ is a smooth curve with isolated singularities, 
having also isolated simple inflections. 
Fix $\xi_0\in I$, such that $\tau_p(\xi_0)\in I_p$ and let $B_P$ be a
Riemannian ball centered at $P:=\phi_p(\xi_0)$ of radius $r$ arbitrary
small.  

The geodesic $\Gamma_{\xi_0}$ divides the boundary of $\partial B_P$ into
two disjoint sectors. 
If $\Sigma_p$ is convex at $P$, the tangent geodesics arriving from
$\gamma$ are entering in the ball $B_P$ crossing one of the two
sectors of the boundary for $0<\xi-\xi_0\ll 1$ and crossing the other
sector for $0<\xi_0-\xi\ll 1$. 
As we have seen in the proof of Proposition \ref{prop:1}, this is
impossible when $\gamma$ has a simple inflection at $\xi_0$. 
Thus $\Sigma_p$ is not convex at $\phi_p(\xi_0)$ if $\gamma$ has a
simple inflection at $\xi_0$. 

The same reasoning allows us to exclude the possibility that
$\Sigma_p$ has a semicubic cusp at $\phi_p(\xi_0)$. 
Hence $\Sigma_p$ has an inflection at  $\phi_p(\xi_0)$. This
inflection is generically simple. 
\end{proof} 

In some degenerate cases the point of $\Sigma_P$ corresponding to a
simple inflection can have degenerate inflections of odd order, like
$t\mapsto(t,t^5)$ in the Euclidian plane, but also more complicated
singularities, like $t\mapsto(t^2,t^4+t^5)$.  

\begin{proposition}\label{prop:6}
The generic singularities of geodesic envelopes of regular curves 
are semicubic cusps, transversal self-intersections and (local) 
second order self-tangencies.  
These singularities are stable under small enough deformations of the
curve $\gamma$ and of the ambien metric (the required smallness of the
deformations depending on the order of the tangential caustics
contening the singularity).  
\end{proposition}
 
The stability concernes only local $2$-self-tangencies: in the
non-local case, two branches of the envelope may have any kind of
tangency, which disappear under small
deformations of the initial curve (splitting into
transversal self-intersections). 

\begin{proof}
This follows from Propositions \ref{prop:4} and \ref{prop:5}. 
\end{proof}

\begin{remark*}
The stability of local $2$-self-tangencies of geodesic envelopes can
also be deduced from the stability of these singularities for
envelopes of tangential families (see \cite{stable}). 
Indeed, denoting by $U$ a connected component of $I_p\smallsetminus
X_p$, the family of geodesics $\{\Gamma_\xi,\xi\in U\}$ is a
tangential family, whose support is the corresponding regular
component of $\Sigma_p$.  
\end{remark*}

We prove now Theorem \ref{thm:2}. 
Consider the manifold $PT^*M$ of all the contact elements on $M$ (a
contact element is a pair $(x,h)$ of a point $x\in M$ and a
homogeneous hyperplane $h\subset T_xM$). 
The {\it projectivized cotangent bundle} $PT^*M$ of $M$ has a natural
contact structure. 

Each regular curve on $M$ can be lifted to a Legendrian curve of
$PT^*M$: at each point of the curve we associate the contact
element formed by the point and the tangent line (on the affine
tangent plane) to the curve at this point.    

By this construction we can associate a surface (called Legendrian
graph, see \cite{legendrian}) 
to the family of the tangent geodesics $\Gamma_\xi$ of our
curve $\gamma$. 
This surface is the union of the Legendrian lifts of
the tangent geodesics.   
The envelope of the family is the (Legendrian) apparent contour of
this Legendrian graph under the natural Legendrian fibration $PT^*M\rightarrow
M$, $(x,h)\mapsto x$.  

Under the hypothesis of Theorem \ref{thm:2}, the contour generator
(i.e. the critical set of this projection) has some regular closed components, 
projecting on the caustics $\Sigma_p$.
Since these components are Legendrian curves, the caustic $\Sigma_p$
can be viewed as a closed fronts. 
A deformation of the initial curve and of the ambient metric induces a
deformation of these Legendrian curves (among Legendrian curves).   
By Proposition \ref{prop:4}, generically (with respect of the curve
$\gamma$) the only singularities of such a front are semicubic cusps
and transversal self-intersections.   
It is proved in \cite{arnold2} that such a front is stable under small
enough deformations of the Legendrian curve generating it.  
Theorem \ref{thm:2} is proved. 

\begin{remark*} 
In the preceding propositions and in Theorem \ref{thm:2} 
the order of the caustic is fixed.  
I do not know whether the same statements hold for all the
tangential caustics simultaneously.  
\end{remark*}



\medskip
\noindent{\sc Gianmarco Capitanio} \\
Department of Mathematical Sciences, Mathematics and Oceanography Building \\
Peach Street,Liverpool. L69 7ZL. United Kingdom \\
{\tt capitani@math.jussieu.fr}\\
{\tt www.institut.math.jussieu.fr/$\sim$capitani}
\end{document}